\documentclass{article}

\usepackage{amsmath,amssymb,amsthm}

\usepackage[latin1]{inputenc}
\usepackage[portuges,english]{babel}
\usepackage{graphicx}
\usepackage{url}

\theoremstyle{definition}
\newtheorem{defn}{Definição}
\theoremstyle{remark}
\newtheorem{exemplo}{Exemplo}

\begin{document}

\selectlanguage{portuges}

\title{Computação Simbólica em Maple\\
no Cálculo das Variações}

% Maple Symbolic Computation in the Calculus of Variations

\author{Andreia M. F. Louro\\
\url{andreia_mfl@yahoo.com.br}
\and Delfim F. M. Torres\\
\url{delfim@ua.pt}}

\date{Centro de Estudos em Optimização e Controlo\\
Departamento de Matemática, Universidade de Aveiro\\
3810-193 Aveiro, Portugal}

\maketitle

% ----------------------------

\begin{abstract}
Neste trabalho pretendemos identificar e ilustrar as potencialidades
e fraquezas do sistema de computação algébrica Maple, na área do
Cálculo das Variações: uma área clássica da Matemática que estuda
os métodos que permitem encontrar valores máximos e mínimos de
funcionais do tipo integral. Os problemas variacionais são
normalmente resolvidos por recurso às condições necessárias de
Euler-Lagrange, que são, em geral, equações diferenciais de segunda
ordem, não lineares e de difícil resolução. Neste trabalho mostramos
como o sistema de computação algébrica Maple pode ser útil na
determinação e resolução das referidas equações. Vamos também
apresentar a resolução do problema clássico da braquistócrona,
sob o ponto de vista do Cálculo das Variações e do Maple, e uma
reformulação do problema restringindo a classe de funções
admissíveis.

\smallskip

\noindent \textbf{Palavras-chave:} sistemas de computação algébrica, Maple,
cálculo das variações, equações de Euler-Lagrange, problema da braquistócrona.
\end{abstract}

\selectlanguage{english}

\begin{abstract}
It is the aim of this work to identify and illustrate
the potential and weaknesses of the computer algebra system Maple
in the area of the Calculus of Variations: a classical area
of mathematics that studies the methods for finding maximum
and minimum values of functionals of integral type. Variational
problems are usually solved with the help of the necessary
optimality conditions of Euler-Lagrange, which are, in general,
nonlinear and difficult second order differential equations
to be solved. We show how the computer algebra system Maple
can be useful in the determination and resolution of these equations.
We will also present the solution to the celebrated brachistochrone problem
from the point of view of the calculus of variations
and the Maple system, and a reformulation
of the classical problem obtained
by restricting the class of admissible functions.

\smallskip


\smallskip

\end{abstract}

\selectlanguage{portuges}

%%%%%%%%%%%%%%%%%%%%%%%%%%%%%%%%%%%%%%%%%%%%%%%%%%%%%%%%%%%%%%%%%%%%%%%%%%%%%%%%%%%%%%%%%%%%%%%%%%%%%%%%%%%%

\section{Introdução}

O \verb"Maple" faz parte de uma família de ambientes computacionais
apelidados de Sistemas de Computação Algébrica. A
computação algébrica, também chamada de computação
simbólica, é uma área de investigação moderna, que surgiu na segunda
metade do século XX. O \verb"Maple" é uma ferramenta matemática
muito poderosa que permite realizar uma miríade de cálculos
simbólicos. Inclui um enorme número de comandos, disponíveis em
vários \emph{packages}, os quais nos permitem trabalhar em áreas
como: Cálculo das Variações, Álgebra Linear, Equações Diferenciais,
Geometria, Lógica, Estatística, entre outras. Inclui também a sua
própria linguagem de programação de alto nível que permite, ao
utilizador, definir os seus próprios comandos.

O Cálculo das Variações é uma área clássica da matemática,
com mais de três séculos de idade e que continua extremamente
activa em pleno século vinte e um. Na sua essência o objectivo
do cálculo das variações é encontrar um caminho, uma curva
ou uma superfície, para os quais uma determinada funcional
tem um valor mínimo ou máximo.

\emph{Nota técnica.} Neste trabalho usámos a versão~9.5 do \verb"Maple".
Ao usar uma versão anterior ou posterior à 9.5, o leitor poderá obter
os resultados numa forma não completamente idêntica à aqui apresentada.
No essencial tudo se manterá, no entanto, igual (verificámos isso mesmo usando,
à posteriori, a versão~11 do \verb"Maple"). À data presente a versão mais recente
do \verb"Maple" é a versão~12 (\textrm{cf.} \url{http://www.maplesoft.com}).

%%%%%%%%%%%%%%%%%%%%%%%%%%%%%%%%%%%%%%%%%%%%%%%%%%%%%%

\section{Cálculo das Variações}

O problema fundamental do cálculo das variações consiste na
determinação de um extremante -- um minimizante ou um maximizante --
para uma funcional integral que depende da escolha de uma função
pertencente a uma determinada classe de funções, em particular
funções cujos valores nos extremos de um dado intervalo real fechado
são fixos.
\\

Sejam $a$, $b$, $a < b$, $A$, $B$ números reais e $L:[a,b] \times
\mathbb{R} \times \mathbb{R} \longrightarrow \mathbb{R}$ uma função
contínua e diferenciável com continuidade. Considera-se o problema seguinte:
\begin{equation}
\label{eqprobfund} J[y(\cdot)] =
\int_{a}^{b}{L\left(x,y(x),y'(x)\right)}dx \longrightarrow \min, \quad
y(a)=A, \quad y(b)=B,
\end{equation}
onde o minimizante se procura no conjunto das funções duas vezes
continuamente diferenciáveis no intervalo $[a,b]$: $y(\cdot) \in
C^{2}([a,b])$. A este problema dá-se o nome de Problema Fundamental
do Cálculo das Variações.
\\

Os problemas deste tipo são normalmente resolvidos recorrendo à
equação diferencial de Euler-Lagrange:
\begin{equation}
\label{eqEL} \frac{\partial L}{\partial
y}(x,y(x),y'(x))-\frac{d}{dx}\frac{\partial L}{\partial
y'}(x,y(x),y'(x))=0 \, .
\end{equation}

\begin{defn}
Às soluções das equações de Euler-Lagrange (\ref{eqEL}) dá-se o nome
de \emph{extremais}.
\end{defn}

De um modo geral, a equação diferencial de Euler-Lagrange é uma
equação não linear, de segunda ordem (ou de ordem superior, quando
os problemas variacionais envolvem derivadas de ordem superior a um)
de difícil resolução.
\\

O package \verb"VariationalCalculus" do
\verb"Maple" proporciona um conjunto de comandos para a resolução de
problemas do Cálculo das Variações. Para utilizar este package temos
de o ``carregar'', executando para isso o seguinte comando:
\begin{verbatim}
  > with(VariationalCalculus):
\end{verbatim}

Neste package temos disponível a função,
\begin{center}
\verb"EulerLagrange(L,x,y(x))"
\end{center}
que devolve um conjunto com as equações de Euler-Lagrange para uma
funcional do tipo (\ref{eqprobfund}). A expressão algébrica
$L$ representa o Lagrangiano,
$x$ a variável independente e $y(x)$ uma expressão algébrica
desconhecida. Esta função do \verb"Maple" permite lidar com
funcionais com mais de uma variável dependente: nesse caso $y(x) \in
\mathbb{R}^n$ é representado, em \verb"Maple", por uma lista.

Vejamos uma aplicação desta função \verb"Maple" no exemplo que se
segue:

\begin{exemplo}
\label{exemplo1} Encontrar a extremal de Euler-Lagrange
$\tilde{y}(\cdot)$ associada à funcional

\begin{equation*}
J[y(.)]=\int_{-1}^{0}(12xy(x)-(y'(x))^2)dx
\end{equation*}
quando sujeita às condições de fronteira
\begin{equation}
\label{condinex1} y(-1)=1 \quad e \quad y(0)=0.
\end{equation}

Começamos por carregar o package em memória e por definir o
Lagrangiano $L$.
\begin{verbatim}
  > with(VariationalCalculus):
\end{verbatim}
\begin{verbatim}
  > L := (12*x*y(x)-diff(y(x),x)^2);
\end{verbatim}
\begin{equation*}
L:=12\,xy \left( x \right) - \left( {\frac {d}{dx}}y \left( x\right)
 \right) ^{2}
\end{equation*}
As respectivas equações de Euler-Lagrange são facilmente obtidas por
intermédio do \verb"Maple".

Utilizando o comando \verb"EulerLagrange" obtemos um conjunto que
contém a equação de Euler-Lagrange:
\begin{verbatim}
  > eqEL := EulerLagrange(L,x,y(x));
\end{verbatim}
\begin{equation*}
eqEL:=\left\{ 12\,x+2\,{\frac {d^{2}}{d{x}^{2}}}y \left( x \right)
 \right\}
\end{equation*}
De notar que, por defeito, quando o segundo membro da equação é $0$,
o \verb"Maple" devolve a equação sob a forma de expressão.

As extremais são obtidas resolvendo esta equação diferencial com as
condições iniciais (\ref{condinex1}) dadas. Recorremos
ao comando \verb"Maple" \verb"dsolve":
\begin{verbatim}
  > dsolve({op(eqEL),y(-1)=1,y(0)=0}, y(x));
\end{verbatim}
\begin{equation*}
y \left( x \right) =-{x}^{3}
\end{equation*}
\end{exemplo}

O comando \verb"op" foi utilizado para extrair o elemento do
conjunto \textit{eqEL}.
\\

Em alguns casos particulares a equação de Euler-Lagrange admite
\emph{primeiros integrais}, isto é, quantidades $\Phi(x,y(x),y'(x))$
que são constantes ao longo de todas as soluções $y(x)$, $x \in
[a,b]$, da equação de Euler-Lagrange, permitindo, em princípio,
reduzir a ordem da mesma e, consequentemente, simplificar o problema.
\\

O \verb"Maple" devolve estes primeiros integrais em conjunto com a
equação de Euler-Lagrange sempre que é capaz de os encontrar.
\\

Em \cite{Louro} consideramos cinco casos de primeiros
integrais. Aqui veremos dois desses casos:
o integral de momento (no exemplo que se segue) e o integral de energia
(aquando a resolução do Problema da Braquistócrona).
\\

Consideremos o seguinte exemplo: determinar as extremais da
funcional

\begin{equation*}
J[y(x)]=\int_{1}^{2}y'(x)(1+x^2y'(x))dx
\end{equation*}
sujeita às condições $y(1)=3$ e $y(2)=5$.
\\

Como $L=y'(x)(1+x^2y'(x))$ não depende explicitamente de
$y(x)$, então a equação de Euler-Lagrange tem a forma
\begin{equation}
\label{eqEL-2c} \frac{d}{dx}\left(\frac{\partial L}{\partial
y'}(x,y'(x))\right)=0
\end{equation}
pelo que admite o \emph{integral de momento}
\begin{equation}
\label{eq:intMomen} \frac{\partial L}{\partial
y'}\left(x,y'(x)\right)=const \, .
\end{equation}
O comando \verb"Maple" \verb"EulerLagrange" devolve neste caso um conjunto
com dois elementos: a equação de Euler-Lagrange propriamente dita
e o integral de momento. Tendo em conta que num conjunto não importa a ordem
de ocorrência dos elementos, não devemos assumir que os conjuntos
devolvidos pelo \verb"Maple" obedeçam ou mantenham determinada
ordem. De maneira a identificarmos de imediato a ordem dos elementos
devolvidos podemos, através do comando
\verb"convert", obter uma lista -- uma
sequência ordenada de elementos entre parêntesis rectos:

\begin{verbatim}
  > eqEL := convert(
     EulerLagrange(diff(y(x),x)*(1+x^2*diff(y(x),x)),x,y(x)),
     list);
\end{verbatim}
\begin{equation*}
eqEL:=\left[-4\,x{\frac {d}{dx}}y \left( x \right) -2\,{x}^{2}{\frac
{d^{2}}{d{x} ^{2}}}y \left( x \right) ,1+2\,{x}^{2}{\frac {d}{dx}}y
\left( x
 \right) = K_{{1}}\right]
\end{equation*}

O \verb"Maple" devolveu uma lista. O primeiro elemento da lista é a
equação de Euler-Lagrange abaixo designada por \verb"eqEL[1]"; o
segundo elemento é o integral de momento, abaixo designado por
\verb"eqEL[2]". Aplicando o comando \verb"dsolve" à equação de
Euler-Lagrange encontramos facilmente a extremal:
\begin{verbatim}
  > dsolve({eqEL[1],y(1)=3,y(2)=5},y(x));
\end{verbatim}
\begin{equation*}
y \left( x \right) =7-\frac{4}{x}
\end{equation*}

Aplicando o comando \verb"dsolve" ao integral de momento, nada é
devolvido.
\begin{verbatim}
  > dsolve({eqEL[2],y(1)=3,y(2)=5},y(x));
\end{verbatim}
O \verb"Maple" fica provavelmente ``baralhado'' a avaliar e a
determinar as constantes. Em vez de introduzir as condições iniciais
no comando \verb"dsolve", vamos avaliá-las mais tarde.
\begin{verbatim}
  > sol := dsolve(eqEL[2],y(x));
\end{verbatim}
\begin{equation*}
sol:=y \left( x \right) =-{\frac
{-\frac{1}{2}+\frac{1}{2}\,K_{{1}}}{x}}+{\it \_C1}
\end{equation*}

Substituindo os valores de $x$ e $y(x)$, respeitantes às condições
iniciais, na solução geral obtida, \textit{sol}, obtemos duas equações com
duas incógnitas.
\begin{verbatim}
  > sol1 := subs({y(x)=3,x=1},sol);
\end{verbatim}
\begin{equation*}
sol_1:=3=\frac{1}{2}-\frac{1}{2}\,K_{{1}}+{\it \_C1}
\end{equation*}
\begin{verbatim}
  > sol2 := subs({y(x)=5,x=2},sol);
\end{verbatim}
\begin{equation*}
sol_2:=5=\frac{1}{4}-\frac{1}{4}\,K_{{1}}+{\it \_C1}
\end{equation*}

Resolvemos o sistema de equações por intermédio do comando
\verb"solve":
\begin{verbatim}
  > const := solve({sol1,sol2},{K[1],_C1});
\end{verbatim}
\begin{equation*}
const:=\left\{ {\it \_C1}=7,K_{{1}}=9 \right\}
\end{equation*}
Finalmente, substituímos o valor das constantes na solução geral,
\textit{sol}, e obtemos a extremal pretendida:
\begin{verbatim}
  > subs({op(const)},sol);
\end{verbatim}
\begin{equation*}
y \left( x \right) =7-\frac{4}{x}
\end{equation*}

Notamos que se derivarmos em ordem a $x$ o integral de momento
obtemos a equação de Euler-Lagrange. Para vermos isto com o
\verb"Maple" utilizamos o comando \verb"diff":
\begin{verbatim}
  > diff(eqEL[2],x);
\end{verbatim}
\begin{equation*}
4\,x{\frac {d}{dx}}y \left( x \right) +2\,{x}^{2}{\frac
{d^{2}}{d{x}^{ 2}}}y \left( x \right) =0
\end{equation*}

Podemos dizer que do ponto de vista matemático a dificuldade de
determinar as extremais por intermédio da equação de Euler-Lagrange
ou do primeiro integral é a mesma, mas do ponto de vista do
\verb"Maple" não: embora equivalentes, o \verb"Maple" não resolve da
mesma forma a equação de Euler-Lagrange e o primeiro integral.
Embora o processo utilizando o integral envolva mais cálculos,
obtivemos, como esperado, a mesma solução.
\\

Para problemas deste tipo, o cálculo das extremais, sob o ponto de
vista do \verb"Maple", revelou-se mais simples recorrendo à equação
de Euler-Lagrange do que utilizando apenas os respectivos primeiros
integrais. No entanto só encontramos uma solução satisfatória para o
problema da Braquistócrona através da resolução a partir do primeiro
integral. O leitor interessado na temática dos primeiros integrais
no cálculo das variações poderá ver \cite{tema}.

%%%%%%%%%%%%%%%%%%%%%%%%%%%%%%%%%%%%%%%%%%%%%%%%%%%%%%

\section{O Problema da Braquistócrona}
\label{braquis}

Um dos primeiros\footnote{O outro problema famoso que esteve
na origem do cálculo das variações e do controlo óptimo
é o problema de Newton da resistência mínima, colocado
e resolvido por Isaac Newton em 1686 \cite{Cristiana3rdJM}.
Para desenvolvimentos recentes deste problema
veja-se \cite{PlakhovTorres,comCristianaCC}.}
e, sem dúvida, o mais famoso problema associado ao
desenvolvimento da teoria matemática do cálculo das variações é,
como já foi referido, o problema da Braquistócrona (para uma
descrição histórica e extremamente interessante deste problema
recomendamos o excelente artigo \cite{Sussmann}):
\\

\textit{Dados dois pontos $A$ e $B$ num plano vertical, qual o
caminho $APB$ que a partícula móvel $P$ atravessa em tempo mínimo,
assumindo que a sua aceleração é apenas devida à gravidade?}
\\

O problema da Braquistócrona é formulado matematicamente como se
segue:
\begin{equation}
\label{braq-mat}
T[y(\cdot)]=\frac{1}{\sqrt{2g}}\int_{x_{0}}^{x_{1}}\sqrt{\frac{1+(y'(x))^2}{y_{0}-y(x)}}dx
\longrightarrow \min
\end{equation}
\begin{equation*}
y(x_{0})=y_{0}, \quad y(x_{1})=y_{1}, \quad y\in C^2(x_{0},x_{1})
\end{equation*}
sendo $T$ o tempo necessário para a partícula deslizar da posição
inicial $A(x_0,y_0)$ até à posição final $B(x_1,y_1)$.
\\

A solução deste problema é um arco de uma ciclóide gerada pelo
movimento de um ponto fixo localizado numa circunferência, de
diâmetro $a$, que rola sobre a recta $y=y_0$ e cujas equações
paramétricas são dadas por
\begin{equation}
\label{eq-param}
\begin{cases}
x=x_{0}+\frac{a}{2}(\theta-\sin(\theta))\,  \\[0.3cm]
y=y_{0}-\frac{a}{2}(1-\cos(\theta)) \,
\end{cases}
\end{equation}
com $\theta_{0}\leq \theta \leq \theta_{1}$, onde $\theta_{0}$ e
$\theta_{1}$ são os valores de $\theta$ nos pontos $A(x_{0},y_{0})$
e $B(x_{1},y_{1})$.
\\

O par de equações

\begin{equation}
\begin{cases}
\label{sist-const}
a\left(\theta_{1}-\sin(\theta_{1})\right)=2(x_{1}-x_{0})\,  \\[0.3cm]
a\left(1-\cos(\theta_{1})\right)=-2(y_{1}-y_{0})
\end{cases}
\end{equation}
permite determinar valores únicos para as constantes $a$ e
$\theta_{1}$ ($0<\theta_{1}<2\pi$) em função dos valores dados
$x_{0}$, $x_{1}$, $y_{0}$ e $y_{1}$.

%%%%%%%%%%%%%%%%%%%%%%%%%%%%%%%%%%%%%%%%%%%%%%%%%%%%%%%%%%%%%%%%%%%%%%%%%%%%%%%%%%%%%%%%%%%%%%%%%%%%%%%

\subsection{Solução da equação de Euler-Lagrange no caso da Braquistócrona}

No Problema da Braquistócrona (\ref{braq-mat}), ignorando o factor
constante $\frac{1}{\sqrt{2g}}$ que não altera em nada o minimizante,
a função $y(\cdot)$ que minimiza
$T$ será encontrada através da equação de Euler-Lagrange
(\ref{eqEL}) com
\begin{equation} \label{braq-lagrangeano}
L=\sqrt{\frac{1+(y'(x))^2}{y_{0}-y(x)}} \, .
\end{equation}
Como $L$ não depende explicitamente de $x$, a equação de
Euler-Lagrange admite o integral de energia

\begin{equation*}
L(y(x),y'(x))-y'(x)\frac{\partial L}{\partial y'}(y(x),y'(x))=const \, ,
\end{equation*}
ou seja,
\begin{equation*}
y'(x)\frac{\partial L}{\partial y'}(y(x),y'(x))-L(y(x),y'(x))=const.
\end{equation*}

Tendo por base ~\cite{Enns} e ~\cite{Smith}, segue-se uma análise do
problema recorrendo aqui ao sistema de computação algébrica
\verb"Maple".

\begin{verbatim}
  > L := (1 + diff(y(x),x)^2)^(1/2)/(y0-y(x))^(1/2);
\end{verbatim}
\begin{equation*}
L:={\frac {\sqrt {1+ \left( {\frac {d}{dx}}y \left( x \right)
\right) ^{ 2}}}{\sqrt {{\it y_0}-y \left( x \right) }}}
\end{equation*}
\begin{verbatim}
  > eqEL := convert(simplify(EulerLagrange(L,x,y(x))),list);
\end{verbatim}
\begin{multline*}
eqEL:=\left[\frac{1}{2}\,{\frac {1+ \left( {\frac {d}{dx}}y \left( x
\right) \right) ^{2 }-2\, \left( {\frac {d^{2}}{d{x}^{2}}}y \left( x
\right) \right) { \it y_0}+2\, \left( {\frac {d^{2}}{d{x}^{2}}}y
\left( x \right)
 \right) y \left( x \right) }{ \left( {\it y_0}-y \left( x \right)
 \right) ^{\frac{3}{2}} \left( 1+ \left( {\frac {d}{dx}}y \left( x \right)
 \right) ^{2} \right) ^{\frac{3}{2}}}}, \right.\\
 \left.{\frac {1}{\sqrt {1+ \left( {\frac {d}{
dx}}y \left( x \right)  \right) ^{2}}\sqrt {{\it y_0}-y \left( x
 \right) }}}=K_{{1}}\right]
\end{multline*}

Neste caso, a equação de Euler-Lagrange possui uma certa
complexidade e o \verb"Maple" não consegue obter uma solução
satisfatória.
\\

Não sendo vantajoso trabalhar com a equação de Euler-Lagrange,
vamos trabalhar com o integral de energia
e determinar com o \verb"Maple" uma forma paramétrica das
extremais para o problema da Braquistócrona.
Ao aplicar o comando \verb"dsolve" ao integral de energia, com a
opção \verb"parametric", o \verb"Maple" devolve as extremais para o
problema na forma paramétrica:
\begin{verbatim}
  > sol := dsolve(eqEL[2],y(x),parametric);
\end{verbatim}
\begin{multline}
\label{sol-maple} sol:=\left[y \left( {\it \_T} \right) ={\frac
{{\it y_0}\,{K_{{1}}}^{2}+{\it y_0} \,{K_{{1}}}^{2}{{\it
\_T}}^{2}-1}{{K_{{1}}}^{2} \left( 1+{{\it \_T}}^{ 2} \right) }} \, , \right.\\
\left.x\left( {\it \_T} \right) ={\frac {{\it \_T}+\arctan
 \left( {\it \_T} \right)
 +\arctan \left( {\it \_T} \right) {{\it \_T}
}^{2}+{\it \_C1}\,{K_{{1}}}^{2}+{\it \_C1}\,{K_{{1}}}^{2}{{\it
\_T}}^{ 2}}{{K_{{1}}}^{2} \left( 1+{{\it \_T}}^{2} \right) }}\right]
\end{multline}

As expressões matemáticas de $x$ e $y$ são dadas em função do
parâmetro ${\it \_T}$. A constante ${\it \_C1}$ é a segunda
constante de integração.
\\

A representação encontrada pelo \verb"Maple" é menos elegante do
que a representação paramétrica usual da ciclóide, atrás descrita.
Conseguimos no entanto obter (\ref{eq-param}), via \verb"Maple", através de alguns artifícios.
\\

Uma nova função $\theta=\theta(x)$ é introduzida através da seguinte relação:
\\

\begin{verbatim}
  > _T := -cot(theta/2);
\end{verbatim}
\begin{equation*}
{\it \_T} := -\cot\left(\frac{1}{2}\,\theta\right)
\end{equation*}
\begin{verbatim}
  > K[1] := 1/sqrt(a);
\end{verbatim}
\begin{equation*}
K_{{1}}:={\frac {1}{\sqrt {a}}}
\end{equation*}

O segundo membro (\emph{right-hand side}) da segunda equação devolvida
em \textit{sol} é dado, em \verb"Maple", por
\verb"rhs(sol[2])"):
\begin{verbatim}
  > x := rhs(sol[2]);
\end{verbatim}
\begin{multline*}
x :=
 \frac{\left( -\cot \left( \frac{\theta}{2} \right) -\frac{\pi}{2} 
 +{\it arccot}\left( \cot \left( \frac{\theta}{2} \right)  \right) 
 + \left( -\frac{\pi}{2} +{
\it arccot} \left( \cot \left( \frac{\theta}{2} \right)  \right)  \right)
 \left( \cot \left( \frac{\theta}{2} \right)  \right) ^{2}\right) a}{1
 + \left( \cot \left( \frac{\theta}{2} \right)  \right) ^{2}} \\
 +\frac{{\it \_C1} + {\it \_C1}\, \left( \cot \left( \frac{\theta}{2}
 \right)  \right) ^{2}}{1
 + \left( \cot \left( \frac{\theta}{2} \right)  \right) ^{2}}
\end{multline*}
A expressão obtida pelo \verb"Maple" é claramente
possível de ser simplificada, bastando, para isso, substituir
${\it arccot}\left(\cot\left(\frac{\theta}{2}\right)\right)$ por $\frac{\theta}{2}$:

\begin{verbatim}
  > x := simplify(subs(arccot(cot(theta/2))=theta/2,x));
\end{verbatim}
\begin{equation*}
x := -\frac{a}{2}\sin\left( \theta \right) +\frac{a}{2}\, \theta-\frac{\pi}{2} \,a+{\it \_C1}
\end{equation*}

Recordando que $\theta=0$ quando $x=x_{0}$, temos ${\it \_C1}=x_{0}+\frac{\pi}{2} a$:
\begin{verbatim}
  > x := subs(_C1=x0+(Pi/2)*a,x);
\end{verbatim}
\begin{equation*}
x :=
-\frac{a}{2}\sin\left(\theta\right) +\frac{a}{2}\,\theta+{\it x0}
\end{equation*}

Seleccionando o segundo membro da primeira equação devolvida em
\textit{sol} (representado, em \verb"Maple", por \verb"rhs(sol[1])")
e aplicando o comando \verb"simplify", obtemos a forma final desejada para $y$:
\begin{verbatim}
  > y := simplify(rhs(sol[1]));
\end{verbatim}
\begin{equation*}
y:=\frac{a}{2}\cos \left( \theta \right)
-\frac{a}{2} + {\it y_0}
\end{equation*}

Com a forma paramétrica de $x$ e $y$ determinada, o tempo de descida mínimo,
$T$, pode ser facilmente calculado em \verb"Maple".
Derivando $x$ e $y$ em ordem a $\theta$, obtemos:

\begin{verbatim}
  > x´(theta) := diff(x,theta);
\end{verbatim}
\begin{equation*}
x'(\theta):=\frac{a}{2} \left( 1-\cos \left( \theta \right)
\right)
\end{equation*}

\begin{verbatim}
  > y´(theta) := diff(y,theta);
\end{verbatim}
\begin{equation*}
y'(\theta):=-\frac{a}{2}\sin \left( \theta \right)
\end{equation*}
De forma a calcular
\begin{equation*}
\frac{\sqrt{1+\left(\frac{\frac{dy}{d\theta}}{\frac{dx}{d\theta}}\right)^2}}{\sqrt{y_{0}-y\left(\theta
\right)}} \, \frac{dx(\theta)}{d\theta}
\end{equation*}
fazemos:
\begin{verbatim}
  > f := simplify(sqrt(1+(y´(theta)/x´(theta))^2)*x´(theta)/sqrt(y0-y));
\end{verbatim}
\begin{equation*}
f:=- \left( -1+\cos \left( \theta \right)  \right) \sqrt {- \left(
-1+ \cos \left( \theta \right)  \right) ^{-1}}{\frac{a}{\sqrt {-a
 \left( -1+\cos \left( \theta \right)  \right) }}}
\end{equation*}

Vamos simplificar a expressão $f$ assumindo $\theta>0$. Designaremos
a expressão devolvida por ``integranda''.
\begin{verbatim}
  > integranda := simplify(f) assuming theta > 0;
\end{verbatim}
\begin{equation*}
integranda:=\sqrt {a}
\end{equation*}

O tempo mínimo de descida, $T$, calcula-se integrando
a ``integranda'' de $\theta=0$ até $\theta_{1}$:
\begin{verbatim}
  > Ttheta1 := (1/sqrt(2*g))*int(integranda,theta=0..theta1);
\end{verbatim}
\begin{equation*}
Ttheta1 := {\frac{\sqrt {2}\sqrt {a}{\it \theta1}}{2\sqrt
{g}}}
\end{equation*}
Os valores das constantes $a$ e $\theta_1$ dependem das
coordenadas dos pontos $A$ e $B$. Para a constante
gravitacional $g$ tomemos o valor aproximado $9.8$:
\begin{verbatim}
  > g := 9.8:
\end{verbatim}

%%%%%%%%%%%%%%%%%%%%%%%%%%%%%%%%%%%%%%%%%%%%%%%%%%%

\subsection{Exemplos}
\label{exe}

\begin{exemplo}
Calculemos o tempo de descida $T$ para o caso em que $A(0,2)$ e
$B(3,1)$ e tracemos a curva de descida que conduz a esse tempo $T$.
\\

Tal como referido anteriormente, o par de equações
(\ref{sist-const}) permite determinar valores únicos para as
constantes $a$ e $\theta_{1}$ ($0<\theta_{1}<2\pi$) em função das
coordenadas $x_{0}$, $y_{0}$ e $x_{1}$, $y_{1}$ dos pontos $A$ e
$B$, respectivamente. Vamos definir funções apropriadas
em \verb"Maple" que nos permitem resolver o problema de forma genérica.
Ilustramos as nossas funções \verb"Maple" para o caso particular
em que $(x_0,y_0) = (0,2)$ e $(x_1,y_1) = (3,1)$.
Determinemos então os valores de $a$ e $\theta_{1}$ para $A(0,2)$ e $B(3,1)$:
\begin{verbatim}
  > eq1 := (x0,x1) -> a*(theta1-sin(theta1))=2*(x1-x0);
\end{verbatim}
\begin{equation*}
      eq1 := (x0, x1) \rightarrow a (\theta1 - \sin(\theta1)) = 2 x1 - 2 x0
\end{equation*}
\begin{verbatim}
  > eq2 := (y0,y1) -> a*(1-cos(theta1))=-2*(y1-y0);
\end{verbatim}
\begin{equation*}
        eq2 := (y0, y1) \rightarrow a (1 - \cos(\theta1))
        = -2 y1 + 2 y0
\end{equation*}
\begin{verbatim}
  > sol := (x0,y0,x1,y1) -> fsolve({eq1(x0,x1),eq2(y0,y1)},{a,theta1});
\end{verbatim}
\begin{equation*}
  sol := (x0, y0, x1, y1) \rightarrow
          fsolve({eq1(x0, x1), eq2(y0, y1)}, {a, \theta1})
\end{equation*}
Os valores de $a$ e $\theta_{1}$ para o nosso exemplo
são obtidos chamando a função \texttt{sol} com os valores
de $x_0 = 0$, $y_0 = 2$, $x_1 = 3$, $y_1 = 1$:
\begin{verbatim}
  > sol(0,2,3,1);
\end{verbatim}
\begin{equation*}
\left\{ a= 1.239374053 \, , {\it \theta1}= 4.051628024  \right\}
\end{equation*}
Para obter os valores de $a$ e $\theta_1$ separadamente,
definimos as funções $va$ (``valor de $a$'') e $vtheta1$
(``valor de $\theta_1$''):
\begin{verbatim}
  > va := (x0,y0,x1,y1) -> eval(a,sol(x0,y0,x1,y1)):
  > va(0,2,3,1);
\end{verbatim}
\begin{equation*}
1.239374053
\end{equation*}
\begin{verbatim}
  > vtheta1 := (x0,y0,x1,y1) -> eval(theta1,sol(x0,y0,x1,y1)):
  > vtheta1(0,2,3,1);
\end{verbatim}
\begin{equation*}
4.051628024
\end{equation*}

Definimos agora em \verb"Maple" a função $Tminimo(x_0,y_0,x_1,y_1)$ que nos
permite calcular o tempo de descida mínimo $T$
e exprimimos a extremal parametricamente através de
$solx(x_0,y_0,x_1,y_1)$ e $soly(x_0,y_0,x_1,y_1)$:
\begin{verbatim}
  > Tminimo := (x0,y0,x1,y1) -> evalf(eval(Ttheta1,sol(x0,y0,x1,y1))):
  > T := Tminimo(0,2,3,1);
\end{verbatim}
\begin{equation*}
T := 1.018832361
\end{equation*}
\begin{verbatim}
  > solx := (xi0,y0,x1,y1) -> eval(x,sol(xi0,y0,x1,y1) union {x0 = xi0}):
  > solx(0,2,3,1);
\end{verbatim}
\begin{equation*}
0.6196870265\,\theta- 0.6196870265\,\sin \left( \theta \right)
\end{equation*}
\begin{verbatim}
  > soly := (x0,yi0,x1,y1) -> subs(sol(x0,yi0,x1,y1) union {y0 = yi0},y):
  > soly(0,2,3,1);
\end{verbatim}
\begin{equation*}
1.380312974+ 0.6196870265\,\cos \left( \theta \right)
\end{equation*}
Ao longo da curva definida parametricamente por
\begin{equation*}
\begin{cases}
x &= 0.6196870265\,\theta- 0.6196870265\,\sin \left( \theta \right) \, , \\
y &= 1.380312974+ 0.6196870265\,\cos \left( \theta \right) \, ,
\end{cases}
\end{equation*}
a partícula demora aproximadamente $1.02$ segundos a deslizar de $A$
até $B$. Este é o tempo mínimo aproximado entre estes dois pontos.
O leitor é convidado a usar as funções \verb"Maple" que acabámos
de definir nos seus próprios exemplos.
\end{exemplo}

Recorrendo ao \verb"Maple", calculamos também facilmente o tempo de descida
da partícula ao longo de outras curvas que unam $A$ a $B$.
\begin{verbatim}
  > L := (y,y0) -> sqrt((1+diff(y(xi),xi)^2)/(y0-y(xi)));
\end{verbatim}
\begin{equation*}
L := (y,y0) \rightarrow \sqrt {{\frac {1+ \left( {\frac {d}{d\xi}}{\it y} \left( \xi \right)
 \right) ^{2}}{{\it y0}-{\it y} \left( \xi \right) }}}
\end{equation*}
\begin{verbatim}
  > Tempo := (y,x0,x1,y0) -> evalf((1/sqrt(2*g))*int(L(y,y0),xi=x0..x1));
\end{verbatim}
\begin{equation*}
Tempo := (y,x0,x1,y0) \rightarrow evalf\left(\frac{1}{\sqrt{2 g}} \int_{x0}^{x1} L(y,y0) d\xi\right)
\end{equation*}

Por exemplo, se determinarmos o tempo de descida, $T_{1}$, ao longo
da recta $y_1$ definida por $A$ e $B$, verificamos que é superior a $T$.

\begin{verbatim}
  > y1 := x -> -(1/3)*x+2;
\end{verbatim}
\begin{equation*}
y1 := x \rightarrow -\frac{x}{3}+2
\end{equation*}
\begin{verbatim}
  > T1 := Tempo(y1,0,3,2);
\end{verbatim}
\begin{equation*}
T1 := 1.428571428
\end{equation*}

Seguindo a proposta feita em \cite{Mello} (embora com outros
valores), podemos também determinar o tempo de descida, $T_{2}$, ao
longo do arco da circunferência com centro (3,6) e que passa pelos
pontos $A(0,2)$ e $B(3,1)$ e verificar que é maior do que $T$. Esta
circunferência tem por equação reduzida: $(x-3)^2+(y-6)^2=25$.

Consideremos a semicircunferência definida por $y_2$:

\begin{verbatim}
  > y2 := x -> 6-(16-x^2+6*x)^(1/2);
\end{verbatim}
\begin{equation*}
y2:= x \rightarrow 6-\sqrt {16-{{\it x}}^{2}+6\,{\it x}}
\end{equation*}
Calculamos facilmente o correspondente
tempo $T_2$ por intermédio da nossa função \verb"Maple":
\begin{verbatim}
  > T2 := Tempo(y2,0,3,2);
\end{verbatim}
\begin{equation*}
T2 := 1.151743820
\end{equation*}

O \verb"Maple" permite criar representações gráficas de várias
curvas, separadamente ou em simultâneo.
Através do comando \verb"display",
disponível no package \verb"plots" do \verb"Maple", conseguimos
obter uma representação gráfica das três curvas acima consideradas:
ciclóide, recta e arco de circunferência. Para podermos utilizar este
comando temos de carregar o package no qual ele se encontra
inserido:

\begin{verbatim}
  > with(plots):
\end{verbatim}

Definindo $p_1$, $p_2$ e $p_3$ como sendo as ``estruturas plot''

\begin{verbatim}
  > p1 := plot([solx(0,2,3,1),soly(0,2,3,1),
                theta=0..vtheta1(0,2,3,1)],
               scaling=constrained, thickness=2,
               labels=["x","y"], color=blue):
\end{verbatim}
\begin{verbatim}
  > p2 := plot(y1(x1),x1=0..3, thickness=2, labels=["x","y"],
               scaling=constrained, color=black):
\end{verbatim}
\begin{verbatim}
  > p3 := plot(y2(x2),x2=0..3, thickness=2, labels=["x","y"],
               scaling=constrained, color=red):
\end{verbatim}
obtemos então uma representação gráfica das três curvas acima
consideradas (ciclóide, recta e arco da circunferência) através
do seguinte comando \texttt{display}:

\begin{verbatim}
  > display(p1,p2,p3);
\end{verbatim}
\begin{center}
\includegraphics[scale=0.55]{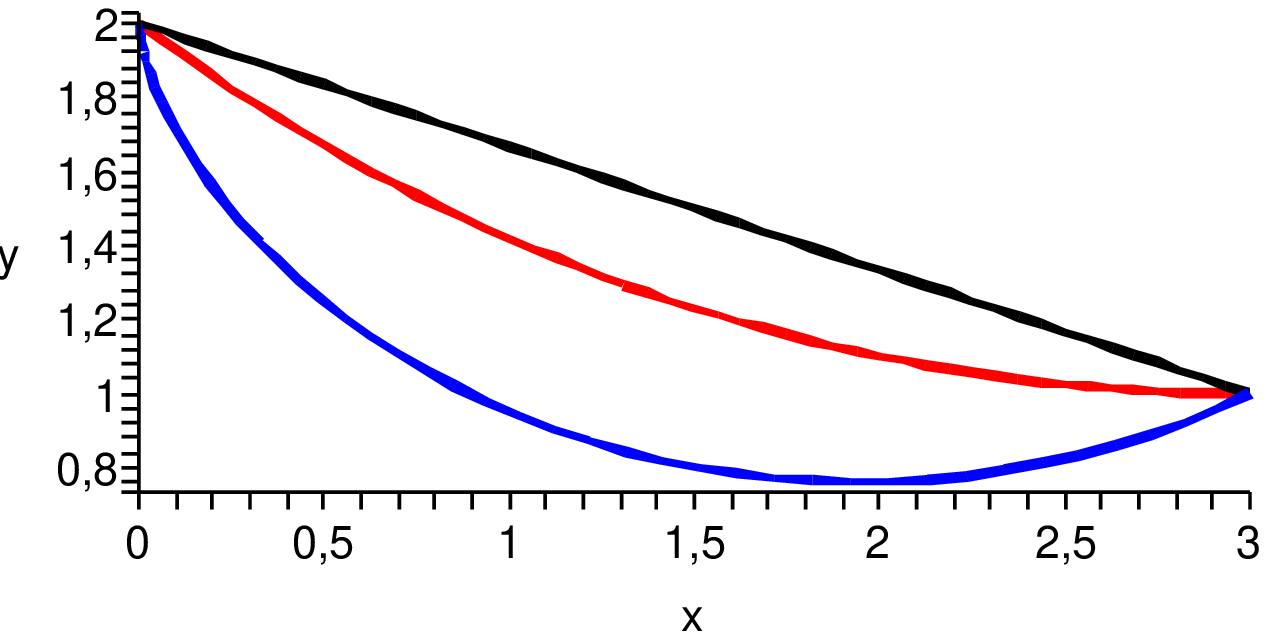}
\end{center}
Embora a ciclóide seja a curva de maior comprimento,
entre as três comparadas, é a curva do tempo mínimo.

\begin{exemplo}
Calculemos o tempo de descida mínimo $T$ para o caso em que $A(1,3)$
e $B(15,1)$ e tracemos a curva.

Para este problema temos $\theta_1$ dado por
\begin{verbatim}
  > vtheta1(1,3,15,1);
\end{verbatim}
\begin{equation*}
4.875635855
\end{equation*}
enquanto $a$ toma o valor de
\begin{verbatim}
  > va(1,3,15,1);
\end{verbatim}
\begin{equation*}
4.776249228
\end{equation*}
O tempo de descida mínimo é inferior a $2.41$ segundos
\begin{verbatim}
  > Tminimo(1,3,15,1);
\end{verbatim}
\begin{equation*}
2.406837209
\end{equation*}
sendo alcançado através da curva descrita parametricamente por
\begin{verbatim}
  > x := solx(1,3,15,1);
\end{verbatim}
\begin{equation*}
x := 1+ 2.388124614\,\theta- 2.388124614\,\sin \left( \theta \right)
\end{equation*}
\begin{verbatim}
  > y := soly(1,3,15,1);
\end{verbatim}
\begin{equation*}
y := 0.611875386+ 2.388124614\,\cos \left( \theta \right)
\end{equation*}

Podemos traçar facilmente esta curva em \verb"Maple" 
usando o comando \texttt{plot} e fazendo variar
o $\theta$ de $0$ até $\theta_1$. Matematicamente, a curva
é parte de uma ciclóide invertida.
\begin{verbatim}
  > plot([x,y,theta=0..vtheta1(1,3,15,1)],
          scaling=constrained, thickness=2, labels=["x","y"]);
\end{verbatim}

\begin{center}
\includegraphics[scale=0.55]{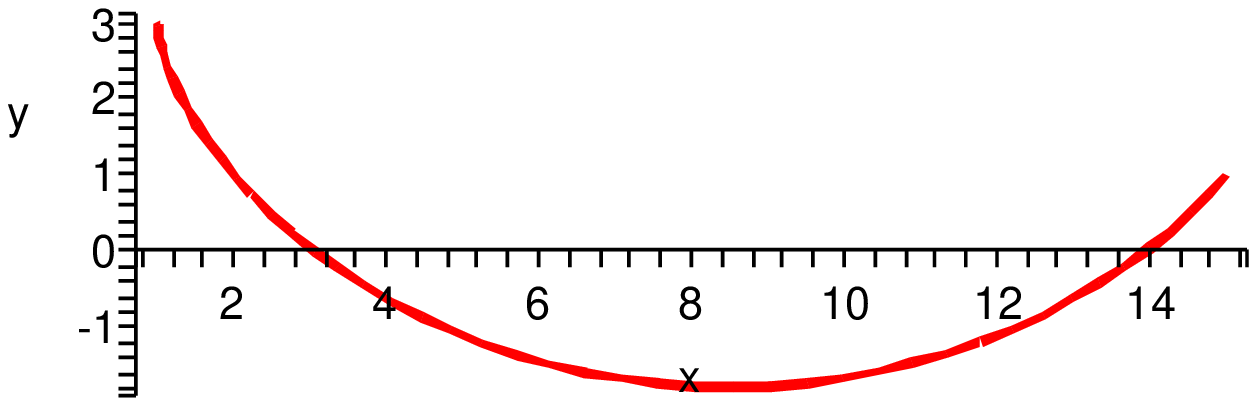}
\end{center}
\end{exemplo}

%%%%%%%%%%%%%%%%%%%%%%%%%%%%%%%%%%%%%%%%%%%%%%%%%%%%%%%%%%%%%%%%%%%%%%%%%%%%%%%

\section{Uma Variante do Problema da Braquistócrona}

Para além do problema publicado por John Bernoulli apresentado em
\S\ref{braquis}, muitas variantes do problema da
Braquistócrona têm sido colocadas e resolvidas ao longo dos tempos:
Braquistócrona através do planeta Terra \cite{Smith};
problema sob a presença de atrito \cite{Haws}; etc.

Na solução clássica foi provado que a Braquistócrona une quaisquer
dois pontos do plano vertical. Isto porque se assumiu a Terra como
inexistente: apenas o campo gravitacional existia.

De acordo com \cite{Ramm}, na reformulação do problema clássico,
que se segue, vamos impor que a linha da Braquistócrona não passe
para valores negativos de $Oy$. Isto implica, obviamente (\textrm{cf.},
por exemplo, a figura no final de \S\ref{exe}),
que a solução geral do problema será diferente.
\\

Sejam $A(0,1)$ e $B(b,0)$, $b>0$. Ignorando o factor constante
$\frac{1}{\sqrt{2g}}$, reconsideramos o problema da Braquistócrona
restringindo a classe das funções admissíveis:
\begin{equation}
\label{braq-mat2}
T[y(\cdot)]=\int_{0}^{1}\sqrt{\frac{1+(y'(x))^2}{1-y(x)}}dx
\longrightarrow \min, \quad y\in S \, ,
\end{equation}
com
\begin{equation*}
S=\left\{y(x): y(0)=1, \quad y(1)=b, \quad y''\geq 0, \quad 0\leq
y(x)\leq y_0(x) \, , \quad y\in C^2(0,b) \right\}
\end{equation*}
onde $y_0(x)$ é a recta que passa pelos pontos $A$ e $B$ e tem por
equação
\begin{equation*}
y_0(x)=1-\frac{x}{b}.
\end{equation*}
Nesta reformulação do problema nem todos os pontos podem ser unidos
por uma extremal da funcional (\ref{braq-mat2}) como se prova em
\cite{Louro} utilizado um novo método de integração da equação de
Euler-Lagrange.

Em ~\cite{Louro}, é explorada também a seguinte questão:

\begin{quotation}
Será verdade que
\begin{equation*}
T(y)\leq T_0:=T(y_0)=2\sqrt{1+b^2}, \quad \forall y \in S \quad?
\end{equation*}
\end{quotation}

Servindo-nos das potencialidades do \verb"Maple" (em particular, recorremos
à função \verb"ThieleInterpolation" disponível no package \verb"CurveFitting"),
encontrámos uma função $\widetilde{y} \in S$
que verifica a desigualdade
\begin{equation*}
T(\widetilde{y})>T_0.
\end{equation*}
O leitor interessado encontrará os detalhes em \cite{Louro}.

%%%%%%%%%%%%%%%%%%%%%%%%%%%%%%%%%%%%%%%%%%%%%%%%%%%%%%%

\section*{Agradecimento}

Agradecemos a um revisor anónimo a leitura cuidada,
as numerosas e pertinentes sugestões de melhoramento
e a experimentação dos comandos aqui apresentados
na versão~10 do \verb"Maple".

%%%%%%%%%%%%%%%%%%%%%%%%%%%%%%%%%%%%%%%%%%%%%%%%%%%%%%%

%%%%%%%%%%%%%%%%%%%%%%%%%%%%%%%%%%%%%%%%%%%%%%%%%%%%%%%


\begin{thebibliography}{9}

\bibitem{Enns} Richard H. Enns,
\emph{Computer Algebra Recipes for Mathematical
Physics}, Birkhäuser, 2005.

\bibitem{tema}
Paulo D. F. Gouveia, Delfim F. M. Torres,
\emph{Computação Algébrica no Cálculo das Variações:
determinação de simetrias e leis de conservação},
TEMA Tend. Mat. Apl. Comput. 6 (2005), no. 1, 81--90.

\bibitem{Haws} LaDawn Haws, Terry Kiser,
\emph{Exploring the Brachistochrone Problem}, The American
Mathematical Monthly, Vol. 102, No. 4, pp. 328-336, 1995.

\bibitem{Louro} Andreia M. F. Louro,
\emph{Computação Simbólica em Maple no Cálculo das Variações},
Dissertação de Mestrado, Univ. Aveiro, 2006.
\url{http://ceoc.mat.ua.pt/dspace/handle/2052/149}

\bibitem{Mello} José L. P. Mello,
\emph{A Rampa de Skate do Tempo Mínimo}, Educação e Matemática --
Revista da Associação de Professores de Matemática, No.~84,
pp.~27--31, 2005.

\bibitem{PlakhovTorres}
Alexander Yu. Plakhov, Delfim F. M. Torres,
\emph{Newton's aerodynamic problem in media
of chaotically moving particles}, Sbornik:
Mathematics, Vol.~196, No.~6, pp.~885--933, 2005.

\bibitem{Ramm} Alexander G. Ramm,
\emph{Inequalities for Brachistochrone}, Mathematical Inequalties
and Applications, Vol.~2, No.~1, pp.~135--140, 1999.

\bibitem{Cristiana3rdJM} Cristiana J. Silva, Delfim F. M. Torres,
\emph{On the Classical Newton's Problem
of Minimal Resistance}, Third Junior
European Meeting on Control, Optimization, and Computation,
University of Aveiro, 6-8 September 2004, Portugal. M.~Guerra and
D.F.M.~Torres eds., Research report CM05/I-04, Dep. Mathematics,
Univ. Aveiro, February 2005, pp.~125--133.

\bibitem{comCristianaCC}  Cristiana J. Silva, Delfim F. M. Torres,
\emph{Two-dimensional Newton's Problem of Minimal Resistance},
Control and Cybernetics, Vol.~35, No.~4, pp.~965--975, 2006.

\bibitem{Smith} Donald R. Smith,
\emph{Variational Methods in Optimization}, Dover Publications
Editora, New York, 1998.

\bibitem{Sussmann} H. J. Sussmann, J. C. Willems,
\emph{300 Anos de Controlo Optimal: da Braquistócrona ao Princípio do Máximo},
Boletim da SPM, 45, pp.~21--54, 2001.

\end{thebibliography}
\end{document}